\documentclass[12pt]{amsart}
\usepackage{amsmath,amsthm,amssymb,amscd}

\begin{document}

\newcommand{\DATE}{November 6, 2002}
\newcommand{\TITLE}{The Rank of Elliptic Surfaces in Unramified Abelian Towers}
\newcommand{\TITLERUNNING}{The Rank of Elliptic Surfaces in Unramified Abelian Towers}
\newcommand{\AUTHOR}{Joseph H. Silverman}

\theoremstyle{plain}
\newtheorem{theorem}{Theorem}
\newtheorem{conjecture}[theorem]{Conjecture}
\newtheorem{proposition}[theorem]{Proposition}
\newtheorem{lemma}[theorem]{Lemma}
\newtheorem{corollary}[theorem]{Corollary}

\theoremstyle{definition}
\newtheorem*{definition}{Definition}

\theoremstyle{remark}
\newtheorem{remark}{Remark}
\newtheorem{example}{Example}
\newtheorem*{acknowledgement}{Acknowledgements}

\newenvironment{notation}[0]{%
  \begin{list}%
    {}%
    {\setlength{\itemindent}{0pt}
     \setlength{\labelwidth}{4\parindent}
     \setlength{\labelsep}{\parindent}
     \setlength{\leftmargin}{5\parindent}
     \setlength{\itemsep}{0pt}
     }%
   }%
  {\end{list}}

\renewcommand{\a}{\alpha}
\renewcommand{\b}{\beta}
\newcommand{\g}{\gamma}
\renewcommand{\d}{\delta}
\newcommand{\e}{\epsilon}
\renewcommand{\l}{\lambda}
\newcommand{\m}{\mu}
\newcommand{\s}{\sigma}
\renewcommand{\t}{\tau}

\newcommand{\F}{\Phi}

\newcommand{\gp}{{\mathfrak{p}}}
\newcommand{\gP}{{\mathfrak{P}}}

\def\Acal{{\mathcal A}}
\def\Bcal{{\mathcal B}}
\def\Ccal{{\mathcal C}}
\def\Dcal{{\mathcal D}}
\def\Ecal{{\mathcal E}}
\def\Fcal{{\mathcal F}}
\def\Gcal{{\mathcal G}}
\def\Hcal{{\mathcal H}}
\def\Ical{{\mathcal I}}
\def\Jcal{{\mathcal J}}
\def\Kcal{{\mathcal K}}
\def\Lcal{{\mathcal L}}
\def\Mcal{{\mathcal M}}
\def\Ncal{{\mathcal N}}
\def\Ocal{{\mathcal O}}
\def\Pcal{{\mathcal P}}
\def\Qcal{{\mathcal Q}}
\def\Rcal{{\mathcal R}}
\def\Scal{{\mathcal S}}
\def\Tcal{{\mathcal T}}
\def\Ucal{{\mathcal U}}
\def\Vcal{{\mathcal V}}
\def\Wcal{{\mathcal W}}
\def\Xcal{{\mathcal X}}
\def\Ycal{{\mathcal Y}}
\def\Zcal{{\mathcal Z}}

\renewcommand{\AA}{\mathbb{A}}
\newcommand{\CC}{\mathbb{C}}
\newcommand{\FF}{\mathbb{F}}
\newcommand{\PP}{\mathbb{P}}
\newcommand{\QQ}{\mathbb{Q}}
\newcommand{\RR}{\mathbb{R}}
\newcommand{\ZZ}{\mathbb{Z}}

\def \bfa{{\mathbf a}}
\def \bfb{{\mathbf b}}
\def \bfc{{\mathbf c}}
\def \bfe{{\mathbf e}}
\def \bff{{\mathbf f}}
\def \bfg{{\mathbf g}}
\def \bfp{{\mathbf p}}
\def \bfr{{\mathbf r}}
\def \bfs{{\mathbf s}}
\def \bft{{\mathbf t}}
\def \bfu{{\mathbf u}}
\def \bfv{{\mathbf v}}
\def \bfw{{\mathbf w}}
\def \bfx{{\mathbf x}}
\def \bfz{{\mathbf z}}


\newcommand{\Aut}{\operatorname{Aut}}
\newcommand{\Cond}{{\mathfrak{N}}} 
\newcommand{\Disc}{\operatorname{Disc}}
\newcommand{\GK}{G_{\Kbar/K}}
\newcommand{\GLK}{G_{L/K}}
\newcommand{\GL}{\operatorname{GL}}
\newcommand{\Image}{\operatorname{Image}}
\newcommand{\Gal}{\operatorname{Gal}}
\newcommand{\Kbar}{{\bar K}}
\newcommand{\MOD}[1]{~(\textrm{mod}~#1)}
\newcommand{\notdivide}{\nmid}
\newcommand{\ord}{\operatorname{ord}}
\newcommand{\rank}{\operatorname{rank}}
\newcommand{\res}{\operatornamewithlimits{res}}
\newcommand{\<}{\langle}
\renewcommand{\>}{\rangle}



\title[\TITLERUNNING]{\TITLE}
\date{\DATE}
\author{\AUTHOR}
\address{Mathematics Department, Box 1917, Brown University,
Providence, RI 02912 USA}
\email{jhs@math.brown.edu}
\thanks{This research supported by NSF DMS-9970382.}

\begin{abstract}
Let $\Ecal\to C$ be an elliptic surface defined over a number
field~$K$.  For a finite covering $C'\to C$ defined over~$K$, let
$\Ecal'=\Ecal\times_CC'$ be the corresponding elliptic surface
over~$C'$. In this paper we give a strong upper bound for the rank
of~$\Ecal'(C'/K)$ in the case of unramified abelien coverings $C'\to
C$ and under the assumption that the Tate conjecture is true
for~$\Ecal'/K$. In the case that~$C$ is an elliptic curve and the map
$C'=C\to C$ is the multiplication-by-$n$ map, the bound
for~$\rank(\Ecal'(C'/K))$ takes the form~$O\left(n^{\kappa/\log\log
n}\right)$, which may be compared with the elementary bound of~$O(n^2)$.
\end{abstract}

\maketitle


\section{Introduction}

It is a longstanding problem to describe the variation of the rank of
the Mordell-Weil group in families of elliptic curves.
There are many variations on this theme. One can study
elliptic curves over number fields or over function fields and one can
study the rank for a fixed base field and varying elliptic curve or
for a fixed elliptic curve and varying base field.
\par
In this note we begin with a number field~$K$, a curve~$C/K$ defined
over~$K$, and a (nonconstant) elliptic
curve~$\Ecal$  defined over the function field~$K(C)$ of~$C$.
Equivalently, we consider an  elliptic surface $\Ecal\to C$
defined~$K$.  Our main result is the following (conditional) upper
bound for the rank of~$\Ecal(K(C'))$ for finite abelian unramified
extensions~$K(C')/K(C)$.

\begin{theorem}
\label{intro:theorem:rankbd}
With notation as above, let $A=\Gal\bigl(K(C')/K(C)\bigr)$. 
Notice that $\GK$ acts on~$A$. Assume that the
Tate conjecture is true for the elliptic surface associated to~$\Ecal/C'$,
and let~$\Cond(\Ecal')$ denote its conductor.
Then
$$
  \rank \Ecal(K(C')) \le \frac{\textup{(Number of $\GK$ orbits of $A$)}}{|A|}
        \cdot
      \left(\bigl|\Cond(\Ecal')\bigr| + 4g' - 4 \right),
$$
where~$g'$ is the genus of the curve~$C'$.
\end{theorem}

We remark that there is an elementary geometric upper bound for the
rank, coming from cohomological considerations, which says that
\begin{equation}
  \label{eqn:geometricrankbound}
  \rank \Ecal(K(C')) \le
      \bigl|\Cond(\Ecal')\bigr| + 4g' - 4 .
\end{equation}
Thus the gain in Theorem~\ref{intro:theorem:rankbd} comes from
nontrivial action of~$\GK$ on the Galois group of the field
extension~$K(C')/K(C)$, or equivalently from nontrivial action of~$\GK$ on
the group of deck transformations of the finite unramified covering
$C'\to C$.
\par
As an interesting special case of the theorem, we take~$C=C'$ to
be an elliptic curve and consider the coverings $[n]:C\to
C$ given by the multiplication-by-$n$ maps. Then $A=C[n]$ is the group
of $n$~torsion points on~$C$ and Serre's theorem tells us that the
action of~$\GK$ on~$A$ is highly nontrivial as~$n$ increases. Combining
this with the theorem gives our second main result.

\begin{theorem}
\label{intro:theorem:ellipticbaserankbound}
Let $C/K$ be an elliptic curve defined over a number field, let $\Ecal/K(C)$
be an elliptic curve, and for each $n\ge1$, let~$K_n$ be the extension field
of~$K(C)$ corresponding to the multiplication-by-$n$ map $[n]:C\to C$. 
Assume that Tate's conjecture is true for the elliptic surface~$\Ecal_n$
associated to~$\Ecal/K_n$. Then there is a~$\kappa>0$ and an
$n_0=n_0(K,C,\d)$ such that
$$
  \rank \Ecal(K_n) \le \bigl|\Cond(\Ecal_n)\bigr|^{\kappa/\log\log n}
  \qquad\text{for all $n>n_0$.}
$$
\end{theorem}

We also show in this situation that the average rank
of~$\Ecal(K_n)$ is smaller than a multiple of the logarithm of its
conductor. (See Theorem~\ref{theorem:ellipticbaserankbound}.) Thus in
an unramified abelian tower over an elliptic base, the rank grows much
more slowly than the elementary geometric bound.

\par
There are many interesting questions one might ask, for example:
\begin{itemize}
\item
Can the rank go to infinity in an unramfied abelian tower? 
\item
What is the best upper bound for the rank in terms of the conductor?
\end{itemize}

In the case that the number field~$K$ is replaced by a finite
field~$\FF_q$, results of Shioda~\cite{Shioda1}, Brumer~\cite{Brumer}
and Ulmer~\cite{Ulmer1} provide a definitive answer. 

\begin{theorem}
\label{intro:theorem:finitefieldbounds}
Let $\Ecal$ be a nonconstant elliptic curve defined over a
function field~$\FF_q(C)$ over a finite field.
\begin{itemize}
\item[\upshape(a)]
The geometric rank of~$\Ecal$ is bounded by
$$
  \rank \Ecal(\bar\FF_q(C)) \le |\Cond(\Ecal)|+4g-4.
$$
Further, there exist examples with conductor of arbitrarily
high degree for which this bound is sharp. (See~\cite{Shioda1,Ulmer1}.)
\item[\upshape(b)]
The arithmetic rank of~$\Ecal$ is bounded by
$$
  \rank \Ecal(\FF_q(C)) \le \frac{|\Cond(\Ecal)|+4g-4}{2\log_q|\Cond(\Ecal)|}
     + O\left(\frac{|\Cond(\Ecal)|}{(\log_q|\Cond(\Ecal)|)^2}  \right).
$$
(See~\cite{Brumer}.)
Further, there exist examples with conductor of arbitrarily
high degree for which the main term in this bound is sharp.
(See~\cite{Ulmer1}.)
\end{itemize}
\end{theorem}

See \cite[Section~4.3]{Ulmer2} for a detailed discussion of these
results.  We also observe that Brumer's proof of the upper
bound (Theorem~\ref{intro:theorem:finitefieldbounds})
is modeled after a result of
Mestre~\cite{Mestre} that deals with elliptic curves over~$\QQ$.

Returning now to the case of an elliptic curve~$\Ecal$ over~$K(C)$
when~$K$ is a number field, we note that the geometric bound given
in~\eqref{eqn:geometricrankbound} holds more generally when the number
field~$K$ is replaced by its algebraic closure, that is,
$$
    \rank \Ecal(\Kbar(C')) \le
      \bigl|\Cond(\Ecal')\bigr| + 4g' - 4 .
$$
This is analogous to the bound in
Theorem~\ref{intro:theorem:finitefieldbounds}(a).
Ulmer~\cite[Section~8]{Ulmer2} has asked if this bound can be improved
for the group~$\Ecal(K(C'))$.  For example, one might be tempted to
make the following conjecture as the analog of the bound in
Theorem~\ref{intro:theorem:finitefieldbounds}(b) and of Mestre's
result~\cite{Mestre}.

\begin{conjecture}
\label{conjecture:rankbound}
Let $K(C)$ be the function field of a curve over a number field, and let
$\Ecal/K(C)$ be a non-constant elliptic curve (i.e., $j(\Ecal)\notin K$).
Then
$$
  \rank\Ecal(K(C)) \ll \frac{|\Cond(\Ecal)|}{\log|\Cond(\Ecal)|},
$$
where the implied constant depends on~$K$ and~$C$. 
\par
More precisely, there is an absolute constant~$\a>0$ so that
$$
  \rank \Ecal(K(C)) \le \a \frac{|\Cond(\Ecal)|+4g-4}{\log|\Cond(\Ecal)|}
      \cdot\log|\Disc(K/\QQ)|.
$$
\end{conjecture}

It is unclear to what extent Theorems~\ref{intro:theorem:rankbd}
and~\ref{intro:theorem:ellipticbaserankbound} should be considered as
providing evidence for this conjecture and to what extent they simply
suggest that the rank over unramified towers is unusually small. 


\begin{acknowledgement}
The author would like to thank Philippe Michel for his help in
understanding the estimates provided by Deligne's work, Doug Ulmer for
encouraging the author to extend the work in~\cite{Silverman1} and
for helpful comments on the first draft of this paper, and
Michael Rosen for numerous mathematical discussions.
\end{acknowledgement}

\section{Setup for a Single Elliptic Surface}

\subsection{Notation}
We set the following notation.

\begin{notation}
\item[$K/\QQ$]
a number field.
\item[$q_\gp$]
the norm of an ideal $\gp$ of $K$.
\item[$C/K$]
a smooth projective curve of genus~$g$.
\item[$\Ecal/K$]
a nonconstant elliptic surface $\Ecal\to C$ defined over~$K$. 
In particular, the assumption that $j(\Ecal)\notin K$ implies that~$\Ecal$
does not split as a product even after base extension of~$C$.
\item[$\Ecal(C/K)$]
the group of sections of~$\Ecal\to C$ defined over~$K$.
\item[$\Cond(\Ecal/C)$]
the conductor of the elliptic surface $\Ecal\to C$. The conductor is
a divisor on~$C$; we denote its degree by~$|\Cond(\Ecal/C)|$.
\end{notation}

\begin{remark}
We fix a finite set of primes~$S$ such that for all $\gp\notin S$, the
elliptic surface $\Ecal\to C$ has good reduction at~$\gp$. That is,
$\Ecal\to C$ has a model over the ring of $S$-integers of~$K$. We
enlarge~$S$ further so that for each prime $\gp\notin S$, the
conductor of $\Ecal/C/\FF_\gp$ is the reduction modulo~$\gp$ of the
conductor of $\Ecal/C/K$. We will write~$\sum_\gp$ to mean the sum
over all primes of~$K$ that are not in~$S$. From time to time, we may 
enlarge the set~$S$. 
\end{remark}

\subsection{A Rank Estimate and a Rank Formula}

An elementary upper bound for the rank of~$\Ecal$ can be obtained
from the cohomology of~$\Ecal(\CC)$. We  call this the 
\textit{geometric bound}, and observe that it automatically provides an
upper bound for the rank of~$\Ecal(C/K)$ over the number field~$K$.

\begin{theorem}[Geometric Rank Bound]
\label{theorem:geometricrankbound}
$$
  \rank \Ecal(C/K) \le \rank \Ecal(C/\Kbar) 
  \le \bigl|\Cond(\Ecal/C)\bigr| + 4g - 4.
$$
\end{theorem}
\begin{proof}
See \cite[Corollary~2]{Shioda}.
\end{proof}

In order to improve the geometric bound, we use a local-global formula
originally proposed by Nagao~\cite{Nagao}. The following notation is
needed in order to state the required result.
\par
For each prime ideal~$\gp$ of~$K$ and each point $R\in C(\FF_\gp)$, 
let~$\Ecal_R$ denote the fiber of~$\Ecal$ over~$R$ and let
$$
  a_\gp(\Ecal_R) = 
  \begin{cases}
      q_\gp + 1 - \bigl|\Ecal_R(\FF_\gp)\bigr| 
             & \text{if $\Ecal_R/\FF_\gp$ is smooth,} \\
      0 & \text{if $\Ecal_R/\FF_\gp$ is singular.} \\
  \end{cases}
$$
If~$\Ecal_R/\FF_\gp$ is smooth, then~$a_\gp(\Ecal_R)$ is the trace of
Frobenius and satisfies the usual Weil bound
$|a_\gp(\Ecal_R)|\le2\sqrt{q_\gp}$. The ``average'' of these values over
the fibers will be denoted
$$
  \Acal_\gp(\Ecal/C) = \frac{1}{q_\gp}\sum_{R\in C(\FF_\gp)} a_\gp(\Ecal_R).
$$
The Weil bound gives $|\Acal_\gp(\Ecal/C)|\le2\sqrt{q_\gp}(1+o(1))$
but a theorem of Deligne (with a further improvement by Michel) gives
a much better estimate.

\begin{theorem}
\label{theorem:delignebound}
$$
  \bigl|\Acal_\gp(\Ecal/C)\bigr| \le \bigl|\Cond(\Ecal/C)\bigr|+4g-4
          + O(1/\sqrt{q_\gp}),
$$
where the implied constant depends on~$\Ecal/C/K$, but is independent
of~$\gp$.
\end{theorem}
\begin{proof}
The weaker upper bound
$$
  \bigl|\Acal_\gp(\Ecal/C)\bigr| \le 2\text{(\# of singular fibers)}+4g-4
$$
follows in a straightforward manner from an equidistribution result of
Deligne~\cite{DeligneWeil2}, see also \cite[(3.6.3)]{KatzGSKSMG}. Each
singular fiber contributes at least one to the conductor, so this
yields
$$
  \bigl|\Acal_\gp(\Ecal/C)\bigr| \le 2\bigl|\Cond(\Ecal/C)\bigr|+4g-4
$$
However, for an elliptic surface, since one knows exactly the
monodromy action around the points of semistable reduction,
Michel~\cite[Section~4]{Michel1} explains how to save one factor of
the semistable part of the conductor, which yields the stated
result. See also Michel~\cite{Michel2} and Fisher~\cite{Fisher} for
similar results and extensions.
\end{proof}

The following analytic version of a conjecture of Nagao~\cite{Nagao}
gives a local-global formula for the rank of~$\Ecal(C/K)$.

\begin{theorem}
\label{theorem:rosensilverman}
Assume that the Tate conjecture is true for the surface $\Ecal/K$. Then
$$
  \rank \Ecal(C/K) 
     = \res_{s=1} \sum_{\gp} -\Acal_\gp(\Ecal/C) \frac{\log q_\gp}{q_\gp^s},
$$
where the sum is over all prime ideals of~$K$ not in~$S$.
\end{theorem}
\begin{proof}
This is proven in \cite[Theorem~1.3]{RosenSilverman}.
\end{proof}

\section{Elliptic Surfaces in Unramified Abelian Towers}

We continue with the notation from the previous section. In
particular,~$C$ is a curve and $\Ecal\to C$ is an elliptic surface,
both defined over a number field~$K$. For any finite cover $C'\to C$, we
obtain a new elliptic surface via pullback,
$$
  \Ecal' = \Ecal\times_C C' \longrightarrow C'.
$$
We consider covers $C'\to C$ satisfying the following conditions: 

\begin{itemize}
\item
The curve~$C'$ and the map $C'\to C$ are defined
over~$K$, so~$\Ecal'/C'$ is likewise defined over~$K$. 
\item
The map $C'\to C$ is a Galois cover with abelian Galois group~$A$.
In other words,~$A$ is an abelian subgroup of $\Aut(C'/\Kbar)$ and
the map $C'\to C$ induces an isomorphism $C'/A\cong C$.
\item
The map $C'\to C$ is unramified.
\end{itemize} 

The conductor of an elliptic surface behaves nicely under an unramified
pullback of the base curve. 

\begin{proposition}
\label{proposition:conductorchange}
Let $\Ecal\to C$ be an elliptic surface and let $f:C'\to C$ be an
unramified map. 
\item[\upshape(a)]
The conductor of the pullback
$\Ecal'=\Ecal\times_CC'$ is given by
$$
  \Cond(\Ecal'/C') = f^*(\Cond(\Ecal/C)).
$$
\item[\upshape(b)]
The canonical divisors on~$C$ and~$C'$ are related by
$$
  \Kcal_{C'} = f^*(\Kcal_C).
$$
\item[\upshape(c)]
In particular,
$$
  |\Cond(\Ecal'/C')|+4g'-4 = 
    \deg(f)\bigl( |\Cond(\Ecal/C)| + 4g-4 \bigr),
$$
where~$g$ and~$g'$ are the genera of~$C$ and~$C'$, repsectively.
\end{proposition}
\begin{proof}
(a) Immediate from the fact that a minimal equation for~$\Ecal$ over the
local ring~$\Ocal_R$ of a point~$R\in C$ remains a minimal Weierstrass
equation over the local ring~$\Ocal_P$ for any point~$P\in f^{-1}(R)$,
since~$\Ocal_P$ is an unramified extension of~$\Ocal_R$.
\par\noindent 
(b) This is \cite[Proposition IV.2.3]{Hartshorne} with trivial ramification
divisor.
\par\noindent
(c) Immediate from~(a),~(b), and
$\deg(\Kcal_C)=2g-2$~\cite[IV.1.3.3]{Hartshorne}.
\end{proof}

Each element of~$A$ is an automorphism $a:C'\to C'$. These
automorphisms need not be defined over~$K$, but the fact that~$C'$ is
defined over~$K$ implies that~$\GK$ acts on~$A$. Since~$A$ is finite,
we can choose a finite extension~$L/K$ so that~$\GLK$ acts on~$A$.
Note that~$\GLK$ need not be abelian.

\begin{example}
If $C$ is an elliptic curve, we can take $C'=C$ and use the multiplication
by $n$~map $[n]:C\to C$. In this case the group~$A$ is the group of
$n$-torsion point~$C[n]$ on~$C$, and the action of~$\GK$ on~$A$ is the
usual Galois action on the $n$-torsion points of an elliptic curve.
\end{example}

\begin{example}
More generally, we can embed~$C$ into its Jacobian variety $C\to J$
and let~$C'$ be the pullback of~$C$ via the multiplication-by~$n$ map
$[n]:J\to J$. Then $C'\to C$ is an abelian unramified cover with group
$A=J[n]$ having the natural~$\GK$ action.
\end{example}

\begin{example}
\label{example:coverscomefromjacobians}
There is a partial converse to the previous example. If $C'\to C$ is
any abelian unramified cover, then there is an isogeny $F:J'\to J$ of
their Jacobians so that $C'=F^{-1}(C)$ and
$A=\ker(F)$. See~\cite[Chapter~VI, Section~12, Corollary to
Proposition~11]{Serre}.
\end{example}

\begin{example}
If we drop the requirement that $C'\to C$ be unramified, then an
interesting case is $C'=C=\PP^1$ with the map $T\to T^n$. This is the
situation studied in~\cite{Silverman1}, where it was shown that the 
rank of $\Ecal'(C'/K)$ is bounded by a generalized divisors-of-$n$ function. If
one further puts various sorts of technical restrictions on the
discriminant of~$\Ecal$, then Shioda,
Stiller~\cite{Stiller}, and Fastenberg~\cite{Fastenberg1,Fastenberg2}
have shown that the rank of~$\Ecal'(C'/K)$ is bounded independently
of~$n$.  The techniques in the present paper can be adapted to handle
coverings~$C'\to C$ with limited ramification, but the resulting
estimates are somewhat complicated, so we have opted to restrict
attention to unramified coverings.
\end{example}

\section{Elementary Results about Groups Acting on Sets}
In this section we recall and prove some elementary estimates that
will be required later.

\begin{lemma}
\label{lemma:numGorbits}
Let~$G$ be a finite group that acts on a finite set~$X$. Then
$$
  \frac{1}{|G|}\sum_{\s\in G} \bigl|\{x\in X : \s(x)=x\}\bigr|
  = \textup{(Number of $G$-orbits of $X$)}.
$$
\end{lemma}
\begin{proof} 
Let $\rho$ be the permutation representation of~$G$ acting on~$X$ and
let~$\chi$ be its character~\cite[Section~1.2]{SerreLRFG}. Then
as in~\cite[Section~2.3]{SerreLRFG},
\begin{align*}
  \text{(Number of }&\text{$G$-orbits of $X$)} \\
  &=\text{(Number of times $\rho$ contains the unit representation)} \\
  &=\frac{1}{|G|}\sum_{\s\in G} \chi(\s). 
\end{align*}
This is the desired formula, since from the definition of the
permutation representation, it is clear that $\chi(\s)$ is equal to
the number of elements of~$X$ that are fixed by~$\s$.
\end{proof}

\begin{lemma}
\label{lemma:numHorbits}
Let~$G$ be a finite group that acts on a finite set~$X$, and let~$H$
be a subgroup of~$G$. Then
$$
  \textup{(Number of $H$-orbits of $X$)}
  \le (G:H)\cdot\textup{(Number of $G$-orbits of $X$)},
$$
with equality if and only if~$H_x=G_x$ for every $x\in X$.
(Here $G_x=\{\s\in G:\s(x)=x\}$ is the stabilizer of~$x$, and similarly for~$H_x$.)
\end{lemma}
\begin{proof}
\begin{align*}
  \text{(Number of $H$-orbits of $X$)}
  &=\sum_{x\in X} \frac{1}{|Hx|} \\
  &=\frac{1}{|H|} \sum_{x\in X} |H_x| \\
  &= (G:H) \sum_{x\in X} \frac{|H_x|}{|G|} \\
  &= (G:H) \sum_{x\in X} \frac{1}{|Gx|}\cdot\frac{|H_x|}{|G_x|} \\
  &\le (G:H) \sum_{x\in X} \frac{1}{|Gx|} \qquad\text{since $H_x\subset G_x$,}\\
  &=(G:H)\cdot\textup{(Number of $G$-orbits of $X$)}.
\end{align*}
This proves the desired inequality. Further, we have equality if and only
if $|H_x|=|G_x|$ for every $x\in X$
\end{proof}

\section{Subgroups of $A=\Aut(C'/C)$ and Intermediate Curves}
Each subgroup $B\subset A$ corresponds to a curve~$C_B$ satisfying 
$$
  C'\to C_B\to C\qquad\text{and}\qquad \Aut(C'/C_B)=B.
$$
Equivalently,~$C_B$ is the quotient curve~$C'/B$. Since~$A$ is
abelian, every subgroup is normal, so the covering $C_B\to C$ is also
Galois with automorphism group naturally isomorphic to~$A/B$.
\par
Recall that we have fixed a finite extension~$L/K$ so that~$\GLK$ acts
on~$A$. All of the curves~$C_B$ are defined over~$L$. If~$\GLK$
normalizes~$B$, that is, if
$$
  \text{$b\in B$ and $\s\in\GLK$} \Longrightarrow \s(b)\in B,
$$
then the curve~$C_B$ is defined over~$K$.
\par
In subsequent sections, we reduce the groups~$B$ and the curves~$C_B$
modulo primes~$\gp$ of~$K$. In order to do this, we observe that~$A$
and its subgroups have a natural structure as group schemes
over~$\Kbar$ (cf. Example~\ref{example:coverscomefromjacobians}), and
in fact as group schemes over the field~$L$. They thus extend to group
schemes over the ring of $S'$-integers of~$L$ for some finite set of
primes~$S'$.  We adjoin to the set~$S$ the primes of~$K$ lying below
the primes of~$S'$.  It then makes sense to reduce not only the
curves~$C_B$, but also the automorphism groups $B=\Aut(C'/C_B)$ and
$A/B=\Aut(C_B/C)$, modulo primes not lying above~$S$.  If~$\GLK$
normalizes~$B$, then~$B$ is a group scheme over~$K$, and we may 
reduce modulo primes of~$K$ not in~$S$.

\section{Abelian Unramified Covers over Finite Fields}
We next reduce the unramified abelian covering $C'\to C$ and the
elliptic surfaces~$\Ecal$ and~$\Ecal'$ modulo a prime ideal~$\gp\notin
S$ of~$K$. Adjoining finitely many primes to our set~$S$ of excluded
primes, we may assume the following conditions:
\begin{itemize}
\item
The reduced curves~$C'/\FF_\gp$ and~$C/\FF_\gp$ are nonsingular.
\item
The map $C'\to C$ over~$\FF_\gp$ is unramifed with abelian Galois
group equal to~$A$. In particular, the map is separable and
$A(\bar\FF_\gp)=A(\CC)$.
\item
The elliptic surfaces~$\Ecal/\FF_\gp$ and~$\Ecal'/\FF_\gp$ are nonsingular.
\item
The conductor~$\Cond(\Ecal/C/\FF_\gp)$ of the reduction of~$\Ecal$
modulo~$\gp$ is equal to the reduction modulo~$\gp$ of the global
conductor $\Cond(\Ecal/C/K)$. More generally, we assume that
this is true for each pullback $\Ecal\times_CC_B$ for each subgroup
$B\subset A$. In particular, it is true for~$\Ecal'$, which corresponds
to taking $B=A$.
\end{itemize}

\begin{proposition}
\label{proposition:imagemodp}
Let $\s=\s_\gp$ be the Frobenius map, so~$\s$ generates the Galois
group of~$\bar\FF_\gp/\FF_\gp$.  Let~$B$ be the subgroup of~$A$
defined by
$$
  B=\{\s(a)\circ a^{-1} : a\in A\}.
$$
\begin{itemize}
\item[\upshape(a)]
The group~$B$ is defined over~$\FF_\gp$, and hence the curve $C_B=C'/B$ is
also defined over~$\FF_\gp$.
\item[\upshape(b)]
The image of $C'(\FF_\gp)$ in~$C(\FF_\gp)$ is the same as the image
of $C_B(\FF_\gp)$ in~$C(\FF_\gp)$.
\item[\upshape(c)] 
Let $R\in C(\FF_\gp)$ be a point in the image of~$C'(\FF_\gp)$. Then
there are exactly $|A(\FF_\gp)|$~points of~$C'(\FF_\gp)$ that map
to~$R$. 
\item[\upshape(d)] 
Let~$R\in C(\FF_\gp)$ be a point in the image of~$C_B(\FF_\gp)$. Then
there are exactly~$(A:B)$ points of~$C_B(\FF_\gp)$ that map to~$R$.
\end{itemize}
\end{proposition}
\begin{proof}
(a) It is clear that $\s(B)=B$, so the fact that~$\s$ generates
$\Gal(\bar\FF_\gp/\FF_\gp)$ implies that~$B$ is defined
over~$\FF_\gp$.  Since~$C'$ is defined over~$\FF_\gp$ by assumption,
this in turn implies that the quotient~$C'/B$ is defined
over~$\FF_\gp$,
\paragraph{(b)}
The composition of maps $C'(\FF_\gp)\to C_B(\FF_\gp)\to C(\FF_\gp)$
shows that the image of~$C'(\FF_\gp)$ is contained in the image
of~$C_B(\FF_\gp)$.  Next let~$Q\in C_B(\FF_\gp)$ and choose any point
$P\in C'(\bar\FF_\gp)$ that maps to~$Q$. The fact that~$Q$ is defined
over~$\FF_\gp$ and that $C_B=C'/B$ means that there is an automorphism
$b\in B$ such that $\s(P)=b(P)$.  By definition of~$B$, there is an
automorphism $a\in A$ such that $b=\s(a)\circ a^{-1}$, and hence
$$
  a^{-1}(P) = \s(a)^{-1}\bigl(\s(P)\bigr) = \s\bigl(a^{-1}(P)\bigr).
$$
Thus $a^{-1}(P)$ is fixed by~$\s$, so it is
in~$C'(\FF_\gp)$. Further,~$a^{-1}(P)$ and~$P$ both have the same
image in~$C(\FF_\gp)$, which in turn is the same as the image of~$Q$
(since~$P\in C'$ maps to~$Q\in C_B$). This proves that given any point
$Q\in C_B(\FF_\gp)$, its image in~$C(\FF_\gp)$ is also the image of a 
point in~$C'(\FF_\gp)$, which gives the other inclusion and completes
the proof of~(b).
\paragraph{(c)}
By assumption, there is at least one point $P\in C'(\FF_\gp)$ that
maps to~$R$. The full inverse image of~$R$ is the orbit
$AR=\{a(R):a\in A\}$. The assumption that $C'\to C$ is unramified
implies that~$AR$ consists of $|A|$~distinct points, or equivalently, 
that only the identity element of~$A$ fixes~$R$. Hence
\begin{align*}
   a(R)\in C'(\FF_\gp)
  &\Longleftrightarrow  \s(a(R))=a(R) \\
  &\Longleftrightarrow  \bigl(a^{-1}\circ \s(a)\bigr)(R) = R \\
  &\Longleftrightarrow  \s(a) = a^{-1} \\
  &\Longleftrightarrow  a \in A(\FF_\gp).
\end{align*}
This completes the proof that there are exactly $|A(\FF_\gp)|$~points
of~$C'(\FF_\gp)$ that map to~$R$. 
\paragraph{(d)}
Applying~(c) to the map~$C_B\to C$, we see that 
$|(A/B)(\FF_\gp)|$  points of~$C_B(\FF_\gp)$ are mapped to~$R$. However,
the definition of~$B$ implies that~$\s$ fixes every point
in~$(A/B)(\bar\FF_\gp)$. To see why this is 
true, let \text{$aB\in (A/B)(\bar\FF_\gp)$}. Then
$$
  \s(aB)=\s(a)B=a\s(a)a^{-1}B=aB,
$$
since $\s(a)a^{-1}\in B$. Hence
$$
  |(A/B)(\FF_\gp)| = |(A/B)(\bar\FF_\gp)| = |A|/|B| = (A:B),
$$
which completes the proof of~(d).
\end{proof}

\begin{proposition}
\label{proposition:Apinequality}
$$
  \bigl|\Acal_\gp(\Ecal')\bigr|
  \le \frac{|A(\FF_\gp)|}{|A|} \cdot 
   \left( \bigl|\Cond(\Ecal'/C')\bigr| + 4g' - 4 \right)
  + O(1/\sqrt{q_\gp}).
$$
\end{proposition}

\begin{remark}
Notice that if~$\Gal(\bar\FF_\gp/\FF_\gp)$ acts trivially on~$A$, then
$A(\FF_\gp)=A(\bar\FF_\gp)$, so the upper bound in
Proposition~\ref{proposition:Apinequality} reduces to the generic
upper bound provided by Theorem~\ref{theorem:delignebound}. However,
when the action is nontrivial, then
Proposition~\ref{proposition:Apinequality} may provide a significant
strengthening of the generic bound.
\end{remark}

\begin{proof}
Let~$B=\{\s(a)\circ a^{-1}:a\in A\}$ be the subgroup of~$A$ described
in Proposition~\ref{proposition:imagemodp}. To ease notation, we let 
$$
   C''=C_B=C'/B \qquad \text{and}\qquad \Ecal''=\Ecal_B=\Ecal\times_CC_B.
$$
Thus we have a commutative diagram
$$
\begin{CD}
  \Ecal'(\FF_\gp) @>>> \Ecal''(\FF_\gp) @>>> \Ecal(\FF_\gp) \\
  @VVV @VVV @VVV \\
  C'(\FF_\gp) @>>> C''(\FF_\gp) @>>> C(\FF_\gp). \\
\end{CD}
$$
Proposition~\ref{proposition:imagemodp} tells us that the image
of~$C'(\FF_\gp)$ in~$C(\FF_\gp)$ is the same as the image of~$C''(\FF_\gp)$
in~$C(\FF_\gp)$ and gives us the multiplicity of each map.
\par
We use the fact that if $P\in C'(\FF_\gp)$ maps to $R\in C(\FF_\gp)$,
then the fiber~$\Ecal'_P$ is isomorphic to~$\Ecal_R$, and similarly if
$Q\in C''(\FF_\gp)$ maps to $R\in C(\FF_\gp)$, then $\Ecal''_Q\cong
\Ecal_R$.  We compute
\begin{align*}
  q_\gp\Acal_\gp(\Ecal')
  &= \sum_{P\in C'(\FF_\gp)} a_\gp(\Ecal'_P)
      \qquad\text{by definition of $\Acal_\gp$,} \\
  &= |A(\FF_\gp)| \hskip-1em
     \sum_{R\in \Image[C'(\FF_\gp)\to C(\FF_\gp)]} \hskip-1em a_\gp(\Ecal_R)
      \qquad\text{from Proposition~\ref{proposition:imagemodp}(c),} \\
  &= |A(\FF_\gp)| \hskip-1em
     \sum_{R\in \Image[C''(\FF_\gp)\to C(\FF_\gp)]} \hskip-1em a_\gp(\Ecal_R)
      \qquad\text{from Proposition~\ref{proposition:imagemodp}(b),} \\
  &= \frac{|A(\FF_\gp)|}{(A:B)}
     \sum_{Q\in C''(\FF_\gp)} a_\gp(\Ecal''_Q)
      \qquad\text{from Proposition~\ref{proposition:imagemodp}(d),} \\
  &= q_\gp \cdot \frac{|A(\FF_\gp)|}{(A:B)} \Acal_\gp(\Ecal'')
\end{align*}
Applying the conductor estimate given in
Theorem~\ref{theorem:delignebound} to~$\Ecal''$ and using the
elementary relation (Proposition~\ref{proposition:conductorchange}c)
between the conductors of~$\Ecal''$ and~$\Ecal'$, we obtain the upper
bound
\begin{align*}
  \bigl|\Acal_\gp(\Ecal')\bigr|
  &=  \frac{|A(\FF_\gp)|}{(A:B)} 
      \cdot\bigl|\Acal_\gp(\Ecal'')\bigr| \\
  &\le  \frac{|A(\FF_\gp)|}{(A:B)} 
       \cdot \left( \bigl|\Cond(\Ecal''/C'')\bigr| + 4g'' - 4 \right) 
       + O(1/\sqrt{q_\gp})\\ 
  &=    \frac{|A(\FF_\gp)|}{(A:B)} \cdot \frac{1}{|B|}
         \cdot \left( \bigl|\Cond(\Ecal'/C')\bigr| + 4g' - 4 \right) 
       + O(1/\sqrt{q_\gp})\\ 
  &=    \frac{|A(\FF_\gp)|}{|A|} 
         \cdot \left( \bigl|\Cond(\Ecal'/C')\bigr| + 4g' - 4 \right)
       + O(1/\sqrt{q_\gp}).
\end{align*}
\end{proof}

\section{An Upper Bound for the Rank}

\begin{theorem}
\label{theorem:orbitupperbound}
Let $\Ecal\to C$ be a be an elliptic surface defined over a number
field~$K$, let $C'\to C$ be an unramified abelian covering defined
over~$K$ and with automorphism group~$A$, and let
$\Ecal'=\Ecal\times_CC'$ be the pullback of~$\Ecal$ via this covering.
Assume that the Tate conjecture is true for~$\Ecal'/K$.
Then
\begin{align*}
  \rank{} &\Ecal'(C'/K) \\
  &\le \frac{\textup{(Number of $\GK$ orbits of $A$)}}
                {|A|}\cdot
      \left(\bigl|\Cond(\Ecal'/C')\bigr| + 4g' - 4 \right).
\end{align*}
\end{theorem}

\begin{remark}
Notice that if~$\GK$ acts trivially on~$A$, then we obtain nothing
better than the geometric upper bound given in
Theorem~\ref{theorem:geometricrankbound}. However, if the Galois
action is nontrivial, as tends to the case in a tower $C\leftarrow
C_1\leftarrow C_2\leftarrow\cdots$, then
Theorem~\ref{theorem:orbitupperbound} often gives an upper bound that
is considerably better than the geometric bound. We will see an
example of this below (Theorem~\ref{theorem:ellipticbaserankbound}) in
which for every $\e>0$, the upper bound for $\rank\Ecal_n(C_n/K)$ is
smaller than $|\Cond(\Ecal_n/C_n)|^\e$ as $n\to\infty$.
\end{remark}

\begin{proof}
We apply the analytic rank formula in Theorem~\ref{theorem:rosensilverman}
to the elliptic surface $\Ecal'\to C'$,
$$
  \rank \Ecal'(C'/K) 
     = \res_{s=1} \sum_{\gp} -\Acal_\gp(\Ecal') \frac{\log q_\gp}{q_\gp^s}.
$$
Taking absolute values and using the estimate for~$|\Acal_\gp|$ provided
by Theorem~\ref{proposition:Apinequality} yields
\begin{align*}
  \rank \Ecal'(C'/K) 
     \le{} &\left(\bigl|\Cond(\Ecal'/C')\bigr| + 4g' - 4 \right)
         \res_{s=1} \sum_{\gp} 
         \frac{|A(\FF_\gp)|}{|A|} \cdot \frac{\log q_\gp}{q_\gp^s} \\
     &+ O\left(\res_{s=1} \sum_{\gp} \frac{\log q_\gp}{q_\gp^{s+1/2}}\right).
\end{align*}
The series in the big-O term converges, so its residue is zero and
it may be discarded. Next we observe that
the size of~$A(\FF_\gp)$ depends only on the action of 
\text{$\gp$-Frobenius} on~$A$. In other words, if we choose an element
$\s\in\GLK$, then~$|A(\FF_\gp)|$ is the same for every
prime~$\gp$ such that~$\s$ is in the \text{$\gp$-Frobenius} conjugacy class
$(\gp,L/K)\subset\GLK$. (As always, we have discarded the finitely many
primes for which~$A$ has bad reduction.) More precisely,
if $\s\in(\gp,L/K)$, then
$$
  |A(\FF_\gp)| = |\{a\in A : \s(a)=a\}|
$$
is simply the number of elements of~$A$ fixed by~$\s$. We denote
this last quantity by $h^0(\s,A)$.
\par
We can thus rewrite the above sum as
$$
  \rank \Ecal'(C'/K) 
     \le \left(\bigl|\Cond(\Ecal'/C')\bigr| + 4g' - 4 \right)
     \sum_{\s\in\GLK} \frac{h^0(\s,A)}{|A|}
      \res_{s=1} \sum_{\gp} \frac{\log q_\gp}{q_\gp^s}.
$$
The residue is the degree of the extension~$L/K$ (cf.~\cite{Silverman1}),
so
$$
  \rank \Ecal'(C'/K) 
     \le \left(\bigl|\Cond(\Ecal'/C')\bigr| + 4g' - 4 \right)
     \frac{1}{|\GLK|}\sum_{\s\in\GLK} \frac{h^0(\s,A)}{|A|}.
$$
Finally, we apply Lemma~\ref{lemma:numGorbits} to obtain the desired
bound, which completes the proof of
Theorem~\ref{theorem:orbitupperbound}.
\end{proof}  

\section{The Rank in an Elliptic Tower}

We now consider the case that the base curve~$C/K$ is an elliptic
curve and we take the unramified abelian covers $[n]:C\to C$ given by
the multiplication-by-$n$ maps. In the notation of
Theorem~\ref{theorem:orbitupperbound}, we have $C'=C$, but~$\Ecal'$ is
most definitely not equal to~$\Ecal$. We let~$\Ecal_n$ denote the
pullback of~$\Ecal$ via the map $[n]:C\to C$. The automorphism
group~$A$ in this case is the group $A=C[n]$ of $n$~torsion points
of~$C$, with the natural action of~$\GK$.  
\par 
As is clear from Theorem~\ref{theorem:orbitupperbound}, the
nontriviality of our bound for the rank of~$\Ecal_n(C/K)$ depends on
the degree of nontriviality of the action of~$\GK$ on~$C[n]$. A famous
theorem of Serre gives us control of that action.

\begin{theorem}[Serre \cite{SerreImageOfGalois}]
\label{theorem:imageofgalois}
Let $C/K$ be an elliptic curve defined over a number field~$K$.
There is an integer~$I(C/K)$ so that for every integer $n\ge1$, the
image of the representation
$$
  \rho_{C,n} : \GK\longrightarrow\Aut(C[n])\cong\GL_2(\ZZ/n\ZZ)
$$
has index at most~$I(C/K)$.
\end{theorem}

\begin{remark}
Conjecturally, the index bound in Theorem~\ref{theorem:imageofgalois}
can be chosen to depend only on the field~$K$, and possibly only on
the degree~$[K:\QQ]$.
\end{remark}

We also need the following elementary fact concerning the
natural action of the general linear group.

\begin{proposition}
\label{proposition:GLorbits}
Let $n\ge1$ and $r\ge1$ be an integers. Then the natural action 
of~$\GL_r(\ZZ/n\ZZ)$ on~$(\ZZ/n\ZZ)^r$ has $d(n)$~distinct orbits,
where~$d(n)$ is the number of divisors of~$n$.
\end{proposition}
\begin{proof}
We begin with the case that $n=p^e$ is a prime power. For any vector
$\bfv=(v_1,\ldots,v_r)\in(\ZZ/p^e\ZZ)^r$, let
$$
  \ord_p(\bfv) = \min\bigl\{\ord_p(v_1),\ldots,\ord_p(v_r),e\bigr\}.
$$
We claim that $\bfv$ and $\bfw$ have the same $\GL_r$-orbit if and only
if $\ord_p(\bfv)=\ord_p(\bfw)$. 
\par
First, if $\bfv=A\bfw$ for any integer matrix~$A$, then it is clear from
the definition that $\ord_p(\bfv)\ge\ord_p(\bfw)$. If~$\bfv$ and~$\bfw$ are in
the same orbit, then~$A$ is invertible, so we get an equality
$\ord_p(\bfv)=\ord_p(\bfw)$. 
\par
Next suppose that $\ord_p(\bfv)=\ord_p(\bfw)$. If this common value
is~$e$, then $\bfv=\bfw=\boldsymbol{0}$ and there is nothing further to
be said. So suppose that the common value is~$k$ with $k<e$. Then we
can write $\bfv=p^k\bfv'$ and $\bfw=p^k\bfw'$, where
$\ord_p(\bfv')=\ord_p(\bfw')=0$.  It thus suffices to prove that any
two vectors with $\ord_p=0$ are in the same $\GL_r(\ZZ/p^e\ZZ)$ orbit,
and for that it suffices to show that if $\ord_p(\bfv)=0$, then~$\bfv$
is in the orbit of the unit vector $\bfe=(1,0,0,\ldots,0)$. We are
thus reduced to showing that every vector with $\ord_p(\bfv)=0$ can be
placed as the first column of a matrix in $\GL_r(\ZZ/p^e\ZZ)$. A
matrix modulo~$p^e$ is invertible if and only if its determinant is
prime to~$p$, so we are reduced to the case that $e=1$.  Then the
condition $\ord_p(\bfv)=0$ says simply the $\bfv\not\equiv{\boldsymbol
0}\MOD{p}$, and the desired conclusion follows from the fact that a
nonzero vector in a vector space can always be extended to a basis for
the vector space. The vectors in this basis, lifted from
$(\ZZ/p\ZZ)^r$ to $(\ZZ/p^e\ZZ)^r$, form the desired matrix
in~$\GL_r(\ZZ/p^e\ZZ)$.
\par
This proves that $\ord_p(\bfv)$ completely characterizes the orbit
of~$\bfv$ under the action of $\GL_r(\ZZ/p^e\ZZ)$. The quantity
$\ord_p(\bfv)$ is an integer between~$0$ and~$e$, which proves that
there are $e+1$ distinct orbits. Since $d(p^e)=e+1$, this proves the
proposition when $n=p^e$ is a prime power. Finally,
the Chinese Remainder Theorem and the multiplicativity of~$d(n)$ 
gives the result for all ingegers $n\ge1$.
\end{proof}

\begin{remark}
Applying Lemma~\ref{lemma:numGorbits} with $G=\GL_r(\ZZ/n\ZZ)$ and
$X=(\ZZ/n\ZZ)^r$ and using Proposition~\ref{proposition:GLorbits}
yields some amusing formulas.  For example, take $r=1$, so
$G=(\ZZ/n\ZZ)^*$, $X=\ZZ/n\ZZ$, and the action is multiplication. Then
for $a\in(\ZZ/n\ZZ)^*$,
\begin{align*}
  \bigl|\{x\in\ZZ/n\ZZ : ax\equiv x \MOD{n} \}\bigr|
  &= \bigl|\{x\in\ZZ/n\ZZ : n|(a-1)x \}\bigr| \\
  &=\gcd(a-1,n)
\end{align*}
This yields
$$
  \sum_{\substack{0\le a<n\\ \gcd(a,n)=1}} \gcd(a-1,n)
  = d(n)\phi(n).
$$
It does not seem obvious, a priori, that there should be any natural
relationship that involves taking numbers that are one less than
numbers relatively prime to~$n$ and looking at their common factors
with~$n$.
\end{remark}

\begin{theorem}
\label{theorem:ellipticbaserankbound}
Let~$C/K$ be an elliptic curve defined over a number field~$K$, 
let $\Ecal\to C$ be an elliptic surface defined over~$K$, and
for each integer $n\ge1$, let $\Ecal_n\to C$ be the elliptic surface
obtained by pullback via the multiplication-by~$n$ map $[n]:C\to C$. 
Assume that the Tate conjecture is true for the surfaces~$\Ecal_n/K$.
\begin{itemize}
\item[\upshape(a)]
Let $I(C/K)$ be the supremum over~$n$ of the index
of $\rho_n(\GK)$ in $\Aut(C[n])$. (Theorem~\ref{theorem:imageofgalois}
ensures that~$I(C/K)$ is finite.) Then for all $n\ge1$,
$$
  \rank \Ecal_n(C/K) \le I(C/K)\cdot \frac{d(n)}{n^2}\cdot |\Cond(\Ecal_n/C)|.
$$
\item[\upshape(b)]
The sum
$$
  \frac{1}{x}\sum_{n\le x} 
      \frac{\rank \Ecal_n(C/K)}{\log\bigl|\Cond(\Ecal_n/C)\bigr|}
$$
is bounded as $x\to\infty$. Thus the average rank of~$\Ecal_n(C/K)$ 
is smaller than a fixed multiple of the logarithm of its conductor.
\item[\upshape(c)]
There is a constant $\kappa=\kappa(K,C,\Ecal)$ so that for all sufficiently
large~$n$,
$$
  \rank \Ecal_n(C/K) 
         \le |\Cond(\Ecal_n/C)|^{\kappa/\log\log|\Cond(\Ecal_n/C)|}.
$$
In particular, for any $\e>0$ we have
$$
  \rank \Ecal_n(C/K) \ll |\Cond(\Ecal_n/C)|^\e,
$$
where the implied constant depends on~$K$,~$C$,~$\Ecal$, and~$\e$, but
is independent of~$n$.
\end{itemize}
\end{theorem}
\begin{proof}
Let
$$
  \rho_n : \GK \longrightarrow \Aut(C[n]) \cong \GL_2(\ZZ/n\ZZ)
$$
be the representation of~$\GK$ on the $n$-torsion of~$C$.  Serre's
theorem (Theorem~\ref{theorem:imageofgalois} tells us that the index
of $\rho_n(\GK)$ in $\GL_2(\ZZ/n\ZZ)$ is at most~$I(C/K)$,
where~$I(C/K)$ is independent of~$n$. It follows from
Proposition~\ref{lemma:numHorbits} and
Proposition~\ref{proposition:GLorbits} that
\begin{align*}
  \text{(Number }&\text{of $\GK$ orbits in $C[n]$)} \\
  &\le I(C/K)\cdot\text{(Number of $\Aut(C[n])$ orbits in $C[n]$)} \\
  &=   I(C/K)\cdot
        \text{(Number of $\GL_2(\ZZ/n\ZZ)$ orbits in $(\ZZ/n\ZZ)^2$)} \\
  &= I(C/K)d(n).
\end{align*}
Applying our main result (Theorem~\ref{theorem:orbitupperbound}) with
$A=C[n]$ and using the above bound for the number of~$\GK$ orbits
in~$C[n]$ yields
$$
  \rank \Ecal_n(C/K) 
     \le I(C/K)\cdot \frac{d(n)}{n^2}\cdot |\Cond(\Ecal_n/C)|.
$$
This completes the proof of~(a). 
\par
Proposition~\ref{proposition:conductorchange} says that
$\Cond(\Ecal_n/C)=n^2\Cond(\Ecal/C)$, so we can rewrite the upper bound
as
$$
  \rank \Ecal_n(C/K) \le I(C/K)\cdot d(n) \cdot |\Cond(\Ecal/C)|.
$$
The divisor function~$d(n)$ has the properties (see \cite[Theorem 3.3
and Theorem 13.12]{Apostol}
$$
  \sum_{n\le x} d(n) \sim x\log x
  \qquad\text{and}\qquad
  \limsup_{n\to\infty} \frac{\log d(n)}{\log n/\log\log n} = \log 2.
$$
The first formula implies that 
$$
  \frac{1}{x} \sum_{2\le n\le x}\frac{d(n)}{\log n}
$$ 
is bounded for all $x\ge2$. Hence 
$$ 
  \frac{1}{x}\sum_{n\le x} 
      \frac{\rank \Ecal_n(C/K)}{\log\bigl|\Cond(\Ecal_n/C)\bigr|}
  \le   \frac{1}{x}\sum_{n\le x} 
      \frac{I(C/K)\cdot d(n) \cdot |\Cond(\Ecal/C)|}
        {\log\bigl(n^2\bigl|\Cond(\Ecal/C)\bigr|\bigr)}
$$
is also bounded.
This completes the proof of~(b). 
Finally, let $c_1,c_2,\ldots$ denote absolute constants. Then
using the second formula gives
\begin{align*}
  \rank \Ecal_n(C/K)
  &\le c_1 I(C/K) |\Cond(\Ecal/C)| n^{c_2/\log\log n} \\
  &\le c_3 I(C/K) |\Cond(\Ecal/C)| 
     |\Cond(\Ecal_n/C)|^{c_4/\log\log |\Cond(\Ecal_n/C)|},
\end{align*}
which completes the proof of~(c).
\end{proof}


\end{document}